\documentclass[12pt]{article}

\usepackage{amsmath,amsthm,amssymb,amsfonts}

\newtheorem{theorem}{Theorem}[section]
\newtheorem{proposition}[theorem]{Proposition}
\newtheorem{corollary}[theorem]{Corollary}
\newtheorem{lemma}[theorem]{Lemma}

\newtheorem{remark}[theorem]{Remark}

\begin{document}

\title{Torus Chiral $n$-Point Functions for Free Boson and Lattice Vertex Operator
Algebras}
\author{Geoffrey Mason\thanks{%
Partial support provided by NSF DMS -9709820 and the Committee on Research,
University of California, Santa Cruz} \\
Department of Mathematics, \\
University of California Santa Cruz, \\
CA 95064, U.S.A. \and Michael P. Tuite\thanks{%
Supported by an Enterprise Ireland Basic Research Grant and the Millenium
Fund, National University of Ireland, Galway } \\
Department of Mathematical Physics, \\
National University of Ireland, \\
Galway, Ireland.\\
and\\
School of Theoretical Physics, \\
Dublin Institute for Advanced Studies, \\
10 Burlington Road, Dublin 4, Ireland.}
\maketitle

\begin{abstract}
We obtain explicit expressions for all genus one chiral $n$-point functions
for free bosonic and lattice vertex operator algebras. We also consider the
elliptic properties of these functions.
\end{abstract}

\section{\protect\smallskip \protect\smallskip Introduction}

\smallskip This is the first of several papers devoted to a detailed and
mathematically rigorous study of chiral $n$-point functions at all genera.
Given a vertex operator algebra (VOA) $V$ (i.e. a chiral conformal field
theory) one may define $n$-point functions at genus one following Zhu \cite
{Z} and use various sewing procedures to define such functions at
successively higher genera. In order to implement such a procedure in
practice, one needs a detailed description of the genus one functions. This
itself is a non-trivial issue, and little seems to be currently rigorously
known beyond certain global descriptions for some specific theories \cite
{DMN}, \cite{DM}. The purpose of the present paper is to supply the needed
information in case $V$ is either a free bosonic or lattice VOA. More
precisely, if $V$ is a free bosonic Heisenberg or lattice VOA, $N$ a $V$%
-module, and $v_{1},\ldots ,v_{n}$ states in $V$, we establish a closed
formula below in Theorem \ref{Big Theorem} for the genus one $n$-point
function $F_{N}(v_{1},z_{1};\ldots ;v_{n},z_{n};\tau )$. Roughly speaking,
in the free boson case the $n$-point functions are elliptic functions whose
detailed structure depends on certain combinatorial data determined by the
states $v_{1},\ldots ,v_{n}$. In the lattice case, the function is naturally
the product of two pieces, one determined by the Heisenberg subalgebra and
one which may be described in terms of the lattice and the genus one prime
form. We note that the role played by elliptic functions and the prime form
in calculating genus one $n$-point functions in string theory has long been
discussed by physicists but a rigorous and complete description of these $n$%
-point functions has been lacking until now e.g. \cite{D,P}.

The paper is organized as follows. We begin in section 2 with a brief review
of relevant aspects of free bosonic Heisenberg and even lattice vertex
operator algebras. Section 3 contains the main results of this work. We
begin with a discussion of free bosonic and lattice VOAs of rank one and
later generalise to the rank $l$ case. We firstly use a recursion formula
for $n$-point functions due to Zhu \cite{Z} to demonstrate that every
lattice $n$-point function is a product of a part determined by the free
bosonic Heisenberg sub-VOA and a part dependent on lattice vectors only. We
also obtain an explicit expression for every free bosonic $n$-point function 
$F_{N}(v_{1},z_{1};\ldots ;v_{n},z_{n};\tau )$ in terms of a combinatoric
sum over specific elliptic functions labelled by data determined by the
states $v_{1},\ldots ,v_{n}$. We next describe the $n$-point functions for
pure lattice states. This involves the identification of such $n$-point
functions as a sum of appropriate weights over a certain set of graphs. This
combinatorial approach then leads to a closed expression for all such $n$%
-point functions in terms of the lattice vectors and the genus one prime
form. Finally we conclude the section with Theorem \ref{Big Theorem} which
describes the expression for every rank $l$ lattice $n$-point function.
Section 4 concludes the paper with a discussion of these $n$-point function
from the point of view of their symmetry and elliptic properties. This
provides some further insight into the nature of the explicit formulas
obtained for $n$-point functions in Section 3.

We collect here notation for some of the more frequently occurring functions
and symbols that will play a role in our work. $\mathbb{N}=\{1,2,3,....\}$ is
the set of positive integers, $\mathbb{Z}$ the integers, $\mathbb{C}$ the complex
numbers, $\mathbb{H}$ the complex upper-half plane. We will always take $\tau $
to lie in $\mathbb{H}$, and $z$ will lie in $\mathbb{C}$ unless otherwise noted.
We set $q_{z}=\exp (z)$ and $q=q_{2\pi i\tau }=\exp (2\pi i\tau )$. For $n$
symbols $z_{1},\ldots ,z_{n}\,$we also set $q_{i}=\exp (z_{i})\,$and $%
z_{ij}=z_{i}-z_{j}$.

We now define some elliptic and modular functions. Let $\wp (z,\omega
_{1},\omega _{2})$ denote the Weierstrass elliptic $\wp $-function with
periods $\omega _{1},\omega _{2}$ and set

\begin{equation}
\wp (z,2\pi i,2\pi i\tau )=\frac{1}{z^{2}}+\sum_{k=4,k\text{ even}}^{\infty
}(k-1)E_{k}(\tau )z^{k-2},  \label{Weierstrass_function}
\end{equation}
so that

\begin{equation}
E_{k}(\tau )=-\frac{B_{k}}{k!}+\frac{2}{(k-1)!}\sum_{n=1}^{\infty }\sigma
_{k-1}(n)q^{n},\quad k\in \mathbb{N},\quad k\text{ even}  \label{Eisensteink}
\end{equation}
is the Eisenstein series of weight $k$ normalized as in \cite{DLM}; $B_{k}$
is a certain Bernoulli number and $\sigma _{k-1}(n)$ a power sum over
positive divisors of $n$. Also set

\begin{equation}
E_{k}(\tau )=0,\quad k\in \mathbb{N},\quad k\text{ odd.}  \label{eisensteinkodd}
\end{equation}
We define

\begin{equation}
P_{0}(z,\tau )=-\log z+\sum_{k=2}^{\infty }\frac{1}{k}E_{k}(\tau )z^{k},
\label{P0}
\end{equation}
related to the genus one prime form $K(z,\tau )$ \cite{Mu} by 
\begin{equation}
K(z,\tau )=\exp (-P_{0}(z,\tau )).  \label{Primeform}
\end{equation}
We further define 
\begin{equation}
P_{n}(z,\tau )=\frac{(-1)^{n}}{(n-1)!}\frac{d^{n}}{dz^{n}}P_{0}(z,\tau )=%
\frac{1}{z^{n}}+\sum_{k=2}^{\infty }\binom{k-1}{n-1}E_{k}(\tau )z^{k-n}
\label{Pndefn}
\end{equation}
Note that $P_{1}(z,\tau )=\varsigma (z,2\pi i,2\pi i\tau )-E_{2}(\tau )z$,
for $\varsigma $ the Weierstrass zeta-function and $P_{2}(z,\tau )=\wp
(z,2\pi i,2\pi i\tau )+E_{2}(\tau )$.

We note two expansions for $P_{2}$:

\begin{eqnarray}
P_{2}(z-w,\tau ) &=&\frac{1}{(z-w)^{2}}+\sum_{r,s\in \mathbb{N}}C(r,s,\tau
)z^{r-1}w^{s-1},  \label{P2expansion1} \\
P_{2}(z+w_{1}-w_{2}) &=&\sum_{r,s\in \mathbb{N}}D(r,s,z,\tau
)w_{1}^{r-1}w_{2}^{s-1}.  \label{P2expansion2}
\end{eqnarray}
(expanding the latter in \smallskip $w_{1},w_{2}$) $\,$so that for $r,s\in 
\mathbb{N}$,

\begin{eqnarray}
C(r,s) &=& C(r,s,\tau )=(-1)^{r+1}\frac{(r+s-1)!}{(r-1)!(s-1)!}E_{r+s}(\tau ),
\label{C(r,s)} \\
D(r,s,z) &=& D(r,s,z,\tau ) =(-1)^{r+1}\frac{(r+s-1)!}{(r-1)!(s-1)!}P_{r+s}(z,\tau ).
\label{D(r,s)}
\end{eqnarray}
We also define for $r\in \mathbb{N}$,

\begin{eqnarray}
C(r,0) &=&C(r,0,\tau )=(-1)^{r+1}E_{r}(\tau ),  \label{C(r,0)} \\
D(r,0,z) &=& D(r,0,z,\tau ) =(-1)^{r+1}P_{r}(z,\tau ).  \label{D(r,0)}
\end{eqnarray}
The Dedekind eta-function is defined by

\begin{equation}
\eta (\tau )=q^{1/24}\prod_{n=1}^{\infty }(1-q^{n}).  \label{etafun}
\end{equation}
Finally, for a (finite) set $\Phi $ we denote by $\Sigma (\Phi )$ the
symmetric group consisting of all permutations of $\Phi $. Set 
\begin{eqnarray}
\mathrm{Inv}(\Phi ) &=&\{\sigma \in \Sigma (\Phi )|\sigma ^{2}=1\},\quad 
\text{(\textit{involutions} of }\Sigma (\Phi )\text{)},\text{ }
\label{InvPhi} \\
\mathrm{Fix}(\sigma ) &=&\{x\in \Phi |\sigma (x)=x\},\quad \text{(\textit{%
fixed-points} of }\sigma \text{)},  \label{Fsigma} \\
F(\Phi ) &=&\{\sigma \in \mathrm{Inv}(\Phi )|\mathrm{Fix}(\sigma )=\emptyset
\},\quad \text{(\textit{fixed-point-free }involutions\textit{\ }of }\Sigma
(\Phi )\text{)}.  \nonumber \\
&&  \label{FPhi}
\end{eqnarray}

\section{Vertex Operator Algebras}

We discuss some aspects of VOA theory to establish context and notation. For
more details, see \cite{FHL}, \cite{FLM}, \cite{Ka}, \cite{MN}.

A vertex operator algebra (VOA) is a quadruple $(V,Y,\mathbf{1},\omega )$
consisting of a $\mathbb{Z}$-graded complex vector space $V=\bigoplus_{n\in 
\mathbb{Z}}V_{n}$ , a linear map $Y:V\rightarrow (\mathrm{End}V)[[z,z^{-1}]]$ ,
and a pair of distinguished vectors (states): the vacuum $\mathbf{1}$ in $%
V_{0}$ , and the conformal vector $\omega $ in $V_{2}$. We adopt
mathematical rather than physical notation for vertex operators, so that for
a state $v$ in $V$, its image under the $Y$ map is denoted

\begin{equation}
Y(v,z)=\sum_{n\in \mathbb{Z}}v(n)z^{-n-1},  \label{Ydefn}
\end{equation}
with component operators (or Fourier modes) $v(n)\in\mathrm{End}V$ and where 
$Y(v,z).\mathbf{1|}_{z=0}=v(-1).\mathbf{1}=v$. We generally take $z$ to be a
formal variable. A concession to physics notation is made concerning the
vertex operator for the conformal vector $\omega $, where we write

\[
Y(w,z)=\sum_{n\in \mathbb{Z}}L(n)z^{-n-2}. 
\]
The modes $L(n)$ close on the Virasoro Lie algebra of central charge $c$:

\[
\lbrack L(m),L(n)]=(m-n)L(m+n)+(m^{3}-m)\frac{c}{12}\delta _{m,-n}. 
\]
We define the homogeneous space of weight $k$ to be $V_{k}=\{v\in
V|L(0)v=kv\}$ where for $v$ in $V_{k}$ we write $wt(v)=k$. Then as an
operator on $V$ we have

\[
v(n):V_{m}\rightarrow V_{m+k-n-1}. 
\]
In particular, the \textit{zero mode} $o(v)=v(wt(v)-1)$ is a linear operator
on each homogeneous space of $V$.

Next we consider some particular VOAs, namely Heisenberg VOAs (or free boson
theories), and lattice VOAs. We consider an $l$-dimensional complex vector
space (i.e., abelian Lie algebra) $\mathfrak{H}$ equipped with a non-degenerate,
symmetric, bilinear form $(,)$ and a distinguished orthonormal basis $%
a_{1},a_{2},...a_{l}$. The corresponding affine Lie algebra is the
Heisenberg Lie algebra $\mathfrak{\hat{H}}=\mathfrak{H}\otimes \mathbb{C}%
[t,t^{-1}]\oplus \mathbb{C}k$ with brackets $[k,\mathfrak{\hat{H}}]=0$ and

\begin{equation}
\lbrack a\otimes t^{m},b\otimes t^{n}]=(a,b)m\delta _{m,-n}k.
\label{Fockbracket}
\end{equation}
Corresponding to an element $\lambda $ in the dual space $\mathfrak{H}^{*}$ we
consider the Fock space defined by the induced (Verma) module

\[
M^{\lambda }=U(\mathfrak{\hat{H}})\otimes _{U(\mathfrak{H}\otimes \mathbb{C}[t]\oplus 
\mathbb{C}k)}\mathbb{C,} 
\]
where $\mathbb{C}$ is the $1$-dimensional space annihilated by $\mathfrak{H}\otimes
t\mathbb{C}[t]$ and on which $k$ acts as the identity and $\mathfrak{H}\otimes
t^{0} $ via the character $\lambda $; $U$ denotes the universal enveloping
algebra. There is a canonical identification of linear spaces

\[
M^{\lambda }=S(\mathfrak{H}\otimes t^{-1}\mathbb{C}[t^{-1}]), 
\]
where $S$ denotes the (graded) symmetric algebra. The Heisenberg VOA $M$
corresponds to the case $\lambda =0$ and the Fock states

\begin{equation}
v=a_{1}(-1)^{e_{1}}.a_{1}(-2)^{e_{2}}....a_{1}(-n)^{e_{n}}.a_{l}(-1)^{f_{1}}.a_{l}(-2)^{f_{2}}...a_{l}(-p)^{f_{p}}.%
\mathbf{1,}  \label{Fockstate}
\end{equation}
for non-negative integers $e_{i},...,f_{j}$ form a\textit{\ }basis of $%
M$. The vacuum $\mathbf{1}$ is canonically identified with the identity of $%
M_{0}=\mathbb{C}$, while the weight 1 subspace $M_{1}$ may be naturally
identified with $\mathfrak{H}$. The vertex operator corresponding to $h$ in $%
\mathfrak{H}$ is given by

\begin{equation}
Y(h,z)=\sum_{n\in \mathbb{Z}}h(n)z^{-n-1},  \label{Yh}
\end{equation}
where $h(n)$ is the usual operator on $M$. $M$ is a simple VOA.

Next we consider the case of lattice VOAs $V_{L}$ associated to a
positive-definite, even lattice $L$ (\cite{B}, \cite{FLM}). Thus $L$ is a
free abelian group of rank $l$, say, equipped with a positive definite,
integral bilinear form $(,):L\otimes L\rightarrow \mathbb{Z}$ such that $%
(\alpha ,\alpha )$ is even for $\alpha \in L$. Let $\mathfrak{H}$ be the space $%
\mathbb{C}\otimes _{\mathbb{Z}}L$ equipped with the $\mathbb{C}$-linear extension of $%
(,)$ to $\mathfrak{H}\otimes \mathfrak{H}$ and let $M$ be the corresponding
Heisenberg VOA. The Fock space of the lattice theory may be described by the
linear space

\begin{equation}
V_{L}=M\otimes \mathbb{C}[L]=\sum_{\alpha \in L}M\otimes e^{\alpha },
\label{VLdefn}
\end{equation}
where $\mathbb{C}[L]$ denotes the group algebra of $L$ with canonical basis $%
e^{\alpha }$, $\alpha \in L$. $M$ may be identified with the subspace $%
M\otimes e^{0}$ of $V_{L}$, in which case $M$ is a subVOA of $V_{L}$ and the
rightmost equation of (\ref{VLdefn}) then displays the decomposition of $%
V_{L}$ into irreducible $M$-modules. We identify $e^{\alpha }$ with the
element $\mathbf{1}\otimes e^{\alpha }$ in $V_{L}$; each of the elements $%
e^{\alpha }$ is a primary state of weight $(\alpha ,\alpha )/2$. The vertex
operator for $h$ in $\mathfrak{H}$ is again represented in the obvious way by (%
\ref{Yh}). The vertex operator for $e^{\alpha }$ is more complicated (loc.
cit.) and is given by

\begin{eqnarray}
Y(e^{\alpha },z) &=&Y_{-}(e^{\alpha },z)Y_{+}(e^{\alpha },z)e^{\alpha
}z^{\alpha },  \nonumber \\
Y_{\pm }(e^{\alpha },z) &=&\exp (\mp \sum_{n>0}\frac{\alpha (\pm n)}{n}%
z^{\mp n}).  \label{Yealpha}
\end{eqnarray}
(The slight inconsistency in notation is more than compensated by its
convenience). The operators $e^{\alpha }\in \mathbb{C}[L]$ have group
commutator 
\begin{equation}
e^{\alpha }e^{\beta }e^{-\alpha }e^{-\beta }=(-1)^{(\alpha ,\beta )},
\label{ealphacomm}
\end{equation}
and $e^{\alpha },z^{\alpha }$ act on any state $u\otimes e^{\beta }\in V_{L}$
as 
\begin{eqnarray}
e^{\alpha }(u\otimes e^{\beta }) &=&\epsilon (\alpha ,\beta )u\otimes
e^{\alpha +\beta },  \label{ealpha} \\
z^{\alpha }(u\otimes e^{\beta }) &=&z^{(\alpha ,\beta )}u\otimes e^{\beta },
\label{zalpha}
\end{eqnarray}
for cocycle $\epsilon (\alpha ,\beta )=\pm 1$. This cocycle can be chosen so
that \cite{FLM} 
\begin{eqnarray}
\epsilon (\alpha ,\beta +\gamma ) &=&\epsilon (\alpha ,\beta )\epsilon
(\alpha ,\gamma ),  \label{cocycleproduct} \\
\epsilon (\alpha ,-\alpha ) &=&\epsilon (\alpha ,\alpha )=1.
\label{cocycleunity}
\end{eqnarray}
$\,$

In the context of his theory of modular-invariance for $n$-point functions
at genus 1, Zhu introduced in \cite{Z} a second VOA $(V,Y[,],\mathbf{1},%
\tilde{\omega})$ associated to a given VOA $(V,Y(,),\mathbf{1},\omega )$.
This will be important in the present paper, and we review some aspects of
the construction here. The underlying Fock space of the second VOA is the
same space $V$ as the first, moreover they share the same vacuum vector $%
\mathbf{1}$ and have the same central charge. The new vertex operators are
defined by a change of co-ordinates \footnote{%
Concerning the co-ordinate change we follow \cite{DLM} rather than \cite{Z}.
The latter has $z$ replaced by $2\pi iz$ in (\ref{Ysquare}). This leads to
minor discrepancies between the notation in \cite{Z} and the present paper
which should be borne in mind.}, namely

\begin{equation}
Y[v,z]=\sum_{n\in \mathbb{Z}}v[n]z^{-n-1}=Y(q_{z}^{L(0)}v,q_{z}-1),
\label{Ysquare}
\end{equation}
while the new conformal vector $\tilde{\omega}$ is defined to be the state $%
\omega -\frac{c}{24}\mathbf{1}$. We set

\begin{equation}
Y[\tilde{\omega},z]=\sum_{n\in \mathbb{Z}}L[n]z^{-n-2}  \label{Ywtilde}
\end{equation}
and write $wt[v]=k$ if $L[0].v=kv$, $V_{[k]}=\{v\in V|wt[v]=k\}$. States
homogeneous with respect to the first degree operator $L(0)$ are \textit{not}
necessarily homogeneous with respect to $L[0]$. On the other hand, it
transpires (cf. \cite{Z}, \cite{DLM}) that the two Virasoro algebras enjoy
the \textit{same} set of primary states. We have $L[-1]=L(0)+L(-1)$, which
leads to the useful relation 
\begin{equation}
o(L[-1]v)=0,  \label{oLminusone}
\end{equation}
for any state $v$. Inasmuch as the co-ordinate change $z\rightarrow
q_{z}=\exp (z)$ maps the complex plane to an infinite cylinder, we sometimes
refer to the VOA as being '\textit{on the sphere}', or '\textit{on the
cylinder}'.

The Heisenberg VOA $M$ is a simple example where there is not too much
difference between being on the sphere or the cylinder. This is basically
because $M$ is generated by its weight 1 states which are primary for both
Virasoro algebras, and because we have $u[1]v=u(1)v=(u,v)\mathbf{1}$ and the
commutator formula

\begin{equation}
\lbrack u[m],v[n]]=m(u,v)\delta _{m,-n},  \label{bosonsq}
\end{equation}
for weight 1 states $u,v\in M$ (cf. \cite{Z}, \cite{DMN} for more details).

\section{\protect\smallskip Torus $n$-point Functions}

In this section we will consider $n$-point functions at genus one. A general
reference is Zhu's paper \cite{Z}. Let $(V,Y,\mathbf{1},\omega )$ be a
vertex operator algebra as discussed in section 2 with $N$ a $V$-module.
Recall (loc. cit.) that for states $v_{1},\ldots v_{n}\in V$ , the $n$-point
function on the torus determined by $N$ is

\begin{equation}
F_{N}(v_{1},z_{1};\ldots
;v_{n},z_{n};q)=Tr_{N}Y(q_{1}^{L(0)}v_{1},q_{1})\ldots
Y(q_{n}^{L(0)}v_{n},q_{n})q^{L(0)-c/24},  \label{npointfunction}
\end{equation}
where $q_{i}=q_{z_{i}}$, $1\leq i\leq n$, for auxiliary variables $%
z_{1},...,z_{n}$. (\ref{npointfunction}) incorporates some cosmetic changes
compared to \cite{Z}: we have adorned our $n$-point functions with an extra
factor $q^{-c/24}$ and omitted a factor of $2\pi i$ from the variables $%
z_{i} $ (cf. footnote to (\ref{Ysquare})). In case $n=1$, (\ref
{npointfunction}) is the usual trace function which we will often denote by

\begin{equation}
Z_{N}(v_{1},\tau )=Tr_{N}o(v_{1})q^{L(0)-c/24}
\end{equation}
where $o(v_{1})$ again denotes the zero mode, now for the vertex operator 
$Y(v_{1},z)$ acting on $N$. Note that this trace is independent of $z_{1}$.
Taking all $v_{i}=\mathbf{1}$ in (\ref{npointfunction}) yields the genus one
partition function for $N$:

\begin{equation}
Z_{N}(\tau )=Tr_{N}q^{L(0)-c/24}=q^{-c/24}\sum_{m\geq 0}\dim N_{m}q^{m}.
\label{ZN}
\end{equation}
where $N_{m}$ is the subspace of $N$ of homogeneous vectors of conformal
weight $m$.

\smallskip The following result, which we use later, holds:

\begin{lemma}
\label{lemma3.1} For states $v_{1},v_{2},\ldots ,v_{n}$ as above we have
\end{lemma}

\begin{eqnarray}
&&F_{N}(v_{1},z_{1};\ldots ;v_{n},z_{n};q)\nonumber\\
&&=Z_{N}(Y[v_{1},z_{1n}].Y[v_{2},z_{2n}]\ldots
Y[v_{n-1},z_{n-1n}].v_{n},\tau )  \label{Fnziminuszn} \\
&&=Z_{N}(Y[v_{1},z_{1}].Y[v_{2},z_{2}]\ldots Y[v_{n},z_{n}].\mathbf{1},\tau )
\label{Fnz1zn}
\end{eqnarray}
where $z_{ij}=z_{i}-z_{j}$.

\textbf{Proof}: Recall notation from section 2 for vertex operator algebras
on the cylinder. Lemma \ref{lemma3.1} is implicit in \cite{Z}, section 4,
especially eqn. (4.4.21). We will give a direct proof in the case $n=2$
based on the associativity of vertex operators. The general case follows in
similar fashion. Associativity tells us (\cite{FHL}, Proposition 3.3.2) that 
\begin{equation}
Tr_{N}Y(v_{1},z_{1})Y(v_{2},z_{2})q^{L(0)}=Tr_{N}Y(Y(v_{1},z_{12})v_{2}),z_{2})q^{L(0)}.
\label{Yassociativity}
\end{equation}
We also have (\cite{FHL}, eqn.(2.6.4))

\begin{equation}
e^{xL(0)}Y(v,y)e^{-xL(0)}=Y(e^{xL(0)}v,e^{x}y).  \label{L0scaling}
\end{equation}
By (\ref{Yassociativity}) and (\ref{L0scaling}) it follows that the
left-hand-side of (\ref{Fnz1zn}) is equal to

\begin{eqnarray*}
&&Tr_{N}Y(Y(q_{1}^{L(0)}v_{1},q_{1}-q_{2})q_{2}^{L(0)}v_{2},q_{2})q^{L(0)-c/24}
\\
&=&Tr_{N}Y(q_{2}^{L(0)}Y(q_{z_{12}}^{L(0)}v_{1},q_{z_{12}}-1).v_{2},q_{2})q^{L(0)-c/24}
\\
&=&Tr_{N}Y(q_{2}^{L(0)}Y[v_{1},z_{12}].v_{2},q_{2})q^{L(0)-c/24} \\
&=&Z_{N}(Y[v_{1},z_{12}].v_{2},\tau ),
\end{eqnarray*}
as required.

$\smallskip $Finally using \cite{FHL}, eqn.(2.3.17) we have in general that

\begin{equation}
e^{xL(-1)}Y(v,y)e^{-xL(-1)}=Y(v,y+x).  \label{Ttranslation}
\end{equation}
Hence (\ref{Fnz1zn}) follows from 
\begin{eqnarray*}
o(Y[v_{1},z_{12}].v_{2}) &=&o(Y[v_{1},z_{12}].e^{-z_{2}L[-1]}.Y[v_{2},z_{2}].%
\mathbf{1}) \\
&=&o(e^{-z_{2}L[-1]}.Y[v_{1},z_{1}].Y[v_{2},z_{2}].\mathbf{1}) \\
&=&o(Y[v_{1},z_{2}].Y[v_{2},z_{2}].\mathbf{1}),
\end{eqnarray*}
and using (\ref{oLminusone}). $\qed $

Recall from the previous section the notation for states and modules for
Heisenberg and lattice vertex operator algebras. We are going to develop
explicit formulas for $n$-point functions in these cases. The final answer
is quite elaborate, so we begin with the rank $1$ case. The general case
will proceed in exactly the same manner. We fix the following notation: $L$
is a rank $l=1$ even lattice with inner product $(,)$, $M$ the corresponding
Heisenberg vertex operator algebra based on the complex space $\mathfrak{H}=\mathbb{%
C}\otimes _{\mathbb{Z}}L$, $a\in \mathfrak{H}$ satisfies $(a,a)=1$, $N=M\otimes
e^{\beta }$ is a simple $M$-module with $\beta \in L$, $h=(\beta ,\beta )/2$
the conformal weight of the highest weight vector of $N$. We will establish
a closed formula for $n$-point expressions of the form

\begin{equation}
F_{N}(v_{1}\otimes e^{\alpha _{1}},z_{1};\ldots ;v_{n}\otimes e^{\alpha
_{n}},z_{n};q),  \label{nptlattice}
\end{equation}
where $\alpha _{1},\ldots ,\alpha _{n}\in L$ and $v_{1},\ldots ,v_{n}$ are
elements in the canonical Fock basis (\ref{Fockstate}) of $M$ on the
cylinder. Thus, $v_{1}=a[-1]^{e_{1}}a[-2]^{e_{2}}\ldots $, etc. Note that
the individual vertex operators $Y(v_{i}\otimes e^{\alpha _{i}},z_{i})$, $%
1\leq i\leq n$ do not generally act on the module $N$, however their
composite does as long as $\alpha _{1}+\ldots +\alpha _{n}=0$. We will
always assume that this is the case. It transpires that (\ref{nptlattice})
factors as 
\begin{equation}
Q_{N}.F_{N}(\mathbf{1}\otimes e^{\alpha _{1}},z_{1};\ldots ;\mathbf{1}%
\otimes e^{\alpha _{n}},z_{n};q),  \label{QNFN}
\end{equation}
where $Q_{N}$ is independent of the $\alpha _{i}$, and our main task will be
to elucidate the structure of the two factors $Q_{N}$ and $F_{N}$. Our
results generalize the calculations in \cite{DMN}, which dealt with the case 
$n=1,\alpha _{1}=0$.

We turn to the precise description of $Q_{N}$. Consider first a Fock state $%
v\in M$ given by

\begin{equation}
v=a[-1]^{e_{1}}\ldots a[-p]^{e_{p}}.\mathbf{1}  \label{vstate}
\end{equation}
where $e_{1},\ldots ,e_{p\text{ }}$are non-negative integers . The state $v$
is determined by a multi-set or \textit{labelled set}, which consists of $%
e_{1}+e_{2}+\ldots +e_{p}$ elements, the first $e_{1}$ of which are labelled 
$1$, the next $e_{2}$ labelled $2$, etc. In this way, each of the states $%
v_{i}$ in (\ref{nptlattice}) is associated with a labelled set $\Phi _{i}$,
and we let $\Phi =\Phi _{1}\cup \ldots \cup \Phi _{n}$ denote the disjoint
union of the $\Phi _{i}$, itself a labelled set. For convenience we often
specify an element of $\Phi $ by its label: the reader should bear in mind
that this expedient can be misleading because there are generally several
distinct elements with the same label.

An element $\varphi \in \mathrm{Inv}(\Phi )$ (cf. (\ref{InvPhi})) considered
as a permutation of $\Phi $, may be represented as a product of cycles, each
of length $1$ or $2$:

\begin{equation}
\varphi =(r_{1}s_{1})\ldots (r_{b}s_{b})(t_{1})\ldots (t_{c}).
\label{phiInv}
\end{equation}
(\ref{phiInv}) tells us that $\Phi =\{r_{1},s_{1},\ldots
,r_{b},s_{b},t_{1},\ldots ,t_{c}\}$, while $\varphi $ exchanges elements
with labels $r_{i}$ and $s_{i}$, and fixes elements with labels $%
t_{1},\ldots ,t_{c}$. Notice that involutions may produce the same
permutation of labels yet correspond to distinct permutations of $\Phi $. We
will always consider such involutions to be distinct, regardless of labels.

Recall the definitions (\ref{C(r,s)}) to (\ref{D(r,0)}). Let $\Xi $ be a
subset of $\Phi $ with $|\Xi |\leq 2$. If $\Xi =\{r,s\}$ has size $2$ with $%
r\in \Phi _{i},s\in \Phi _{j}$, we define 
\begin{equation}
\gamma (\Xi )=\left\{ 
\begin{array}{c}
D(r,s,z_{ij},\tau ),\quad i\neq j \\ 
C(r,s,\tau ),\quad i=j.
\end{array}
\right.  \label{GammaX}
\end{equation}
Note that $D(r,s,z,\tau )=D(s,r,-z,\tau )$, so that the order in which the
arguments $r,s$ appear is of no relevance. If $\Xi =\{r\}\subseteq \Phi _{k}$
we define

\begin{equation}
\gamma (\Xi )=(a,\delta _{r,1}\beta +C(r,0,\tau )\alpha
_{k}+\sum_{l>k}D(r,0,z_{kl},\tau )\alpha _{l}).  \label{gammaPhik}
\end{equation}
For $\varphi \in \mathrm{Inv}(\Phi )$ set

\begin{equation}
\Gamma (\varphi )=\prod_{\Xi }\gamma (\Xi )  \label{GammaphiInv}
\end{equation}
where the product ranges over all orbits (cycles) of $\varphi $ in its
action on $\Phi $. Finally, set

\begin{equation}
Q_{N}(v_{1},z_{1};\ldots ;v_{n},z_{n};q)=\sum_{\varphi \in \mathrm{Inv}(\Phi
)}\Gamma (\varphi )  \label{Qnv1vn}
\end{equation}
We can now formally state our first main result about $n$-point functions.

\begin{proposition}
\label{Propgennpt} Let $v_{1},\ldots ,v_{n}$ be states of the form (\ref
{vstate}) in the rank $1$ free boson theory $M$ and let $\Phi $ be the
labelled set determined (as above) by these states. Then the following holds
for lattice elements $\alpha _{1},\ldots ,\alpha _{n}$ satisfying $\alpha
_{1}+\ldots +\alpha _{n}$ $=0$:
\end{proposition}

\begin{eqnarray}
&&F_{N}(v_{1}\otimes e^{\alpha _{1}},z_{1};\ldots ;v_{n}\otimes e^{\alpha
_{n}},z_{n};q)  \nonumber \\
&=&Q_{N}(v_{1},z_{1};\ldots ;v_{n},z_{n};q)F_{N}(\mathbf{1}\otimes e^{\alpha
_{1}},z_{1};\ldots ;\mathbf{1}\otimes e^{\alpha _{n}},z_{n};q).
\label{Fnlattice}
\end{eqnarray}

\textbf{Proof:} The idea is to carefully examine a recursion formula for $n$%
-point functions due to Zhu (\cite{Z}, Proposition 4.3.3). Bearing in mind
the differences in notation between the present paper and \cite{Z}, we quote
the following:

\begin{lemma}[Zhu]
\label{LemZhurec} Assume that $u_{1},\ldots ,u_{n},b$ are states in a vertex
operator algebra $V$, that $N$ is a $V$-module, and that $o(b)$ acts as a
scalar on $N$. Then 
\begin{eqnarray}
&&F_{N}(b[-1]u_{1},z_{1};\ldots ;u_{n},z_{n};q)  \nonumber \\
&&=Tr_{N}o(b)Y(q_{1}^{L(0)}u_{1},q_{1})\ldots
Y(q_{n}^{L(0)}u_{n},q_{n})q^{L(0)-c/24}  \nonumber \\
&&+\sum_{k\geq 1}E_{2k}(\tau )F_{N}(b[2k-1]u_{1},z_{1};u_{2},z_{2};\ldots
;u_{n},z_{n};q)  \nonumber \\
&&+\sum_{m\geq 0}\sum_{k=2}^{n}(-1)^{m+1}P_{m+1}(z_{k1},\tau
)F_{N}(u_{1},z_{1};\ldots ;b[m]u_{k},z_{k};...;u_{n},z_{n};q)  \nonumber \\
&&-\frac{1}{2}%
\sum_{k=1}^{n}F_{N}(u_{1},z_{1};...;b[0]u_{k},z_{k};...;u_{n},z_{n};q).
\label{FNb1Zhu}
\end{eqnarray}
\end{lemma}

$\qed $

We apply this result in the case that $N$ and $u_{i}=v_{i}\otimes e^{\alpha
_{i}}$ are as discussed in (\ref{nptlattice}) - (\ref{Qnv1vn}), with $b=a$.
The zero mode of $a$ acts on $N$ as multiplication by the scalar $(a,\beta )$
and $a[0].v_{i}\otimes e^{\alpha _{i}}=(a,\alpha _{i})v_{i}\otimes e^{\alpha
_{i}}$. As a result of $\alpha _{1}+\ldots +\alpha _{n}$ $=0$ it follows
that the last summand in (\ref{FNb1Zhu}) vanishes, and (\ref{FNb1Zhu}) then
reads

\begin{eqnarray}
&&F_{N}(a[-1]v_{1}\otimes e^{\alpha _{1}},z_{1};\ldots ;v_{n}\otimes
e^{\alpha _{n}},z_{n};q)  \nonumber \\
&&=(a,\beta +\sum_{k=2}^{n}D(1,0,z_{1k})\alpha _{k})F_{N}(v_{1}\otimes
e^{\alpha _{1}},z_{1};\ldots ;v_{n}\otimes e^{\alpha _{n}},z_{n};q) 
\nonumber \\
&&+\sum_{k\geq 1}C(2k-1,1)F_{N}(\hat{v}_{1}\otimes e^{\alpha
_{1}},z_{1};\ldots ;v_{n}\otimes e^{\alpha _{n}},z_{n};q)  \nonumber \\
&&+\sum_{m\geq 1}\sum_{k=2}^{n}
D(m,1,z_{k1})F_{N}(v_{1}\otimes
e^{\alpha _{1}},z_{1};\ldots \hat{v}_{k}\otimes e^{\alpha _{k}},z_{k};\ldots
;v_{n}\otimes e^{\alpha _{n}},z_{n};q).  \nonumber \\
&&  \label{FNa1Zhu}
\end{eqnarray}
Here, we have used the (admittedly uninformative) notation $\hat{v}_{1}$ in
the second summand to indicate that a factor $a[-2k+1]$ should be removed
from the expression of $v_{1}$ as a product (\ref{vstate}), and indeed that
this should be implemented as often as $a[-2k+1]$ occurs in the expression.
If $a[-2k+1]$ does not occur in the expression for $v_{1}$ then $\hat{v}_{1}$
is defined to be zero. Similar notation $\hat{v}_{k}$ occurs in the third
summand, where it indicates removal of a factor $a[-m]$.

Next we develop the analog of (\ref{FNa1Zhu}) in which $a[-1]$ is replaced
by $a[-p]$ for any positive integer $p$. To this end we take $%
b=L[-1]^{p-1}.a $. We easily calculate that $b[m]=(-1)^{p-1}\binom{m}{p-1}%
a[m-p+1]$, in particular $b[-1]=a[-p]$ and $b[0]=0$ if $p\geq 2$. Note also
that $o(b)=0$ if $p\geq 2$, thanks to (\ref{oLminusone}). With this choice
of $b$ and $p$, and after some calculation, (\ref{FNb1Zhu}) reduces to the
next equation. In fact, we can combine the resulting equality (for $p\geq 2$%
) with the case $p=1$. What obtains is the basic recursive relation
satisfied by our $n$-point functions, namely 
\begin{eqnarray}
&&F_{N}(a[-p]v_{1}\otimes e^{\alpha _{1}},z_{1};\ldots ;v_{n}\otimes
e^{\alpha _{n}},z_{n};q)  \nonumber \\
&&=\sum_{k>p/2}C(2k-p,p )F_{N}(\hat{v}_{1}\otimes e^{\alpha
_{1}},z_{1};\ldots ;v_{n}\otimes e^{\alpha _{n}},z_{n};q)  \nonumber \\
&&+\sum_{m>p-1}\sum_{k=2}^{n}D(m-p+1,p,z_{k1} ).  \nonumber \\
&&F_{N}(v_{1}\otimes e^{\alpha _{1}},z_{1};\ldots ;\hat{v}_{k}\otimes
e^{\alpha _{k}},z_{k};\ldots ;v_{n}\otimes e^{\alpha _{n}},z_{n};q) 
\nonumber \\
&&+(a,\delta _{p,1}\beta +C(p,0 )\alpha
_{1}+\sum_{k=2}^{n}D(p,0,z_{1k} )\alpha _{k}).  \nonumber \\
&&F_{N}(v_{1}\otimes e^{\alpha _{1}},z_{1};\ldots ;v_{n}\otimes e^{\alpha
_{n}},z_{n};q),  \label{latticerecursion}
\end{eqnarray}
where we have used a similar convention to the case $p=1$ regarding symbols $%
\hat{v}_{1},\hat{v}_{k}$.

Close scrutiny of relation (\ref{latticerecursion}) reveals how to complete
the proof of Proposition \ref{Propgennpt}, which at this point is a matter
of interpreting the recursion formula. We choose an element with label $p$
from the first labelled set $\Phi _{1}$ determined by $v_{1}$. The first sum
on the r.h.s of (\ref{latticerecursion}) then corresponds to certain terms
in the representation (\ref{vstate}) of $v_{1}$. Indeed, as long as $2k>p$,
a factor $a[p-2k]$ will give rise to a term $C(2k-p,p,\tau )F_{N}(...)$, and
via (\ref{GammaX}) we identify $C(2k-p,p,\tau )$ with $\gamma (\Xi )$ where $%
\Xi =\{2k-p,p\}$ $\subseteq $ $\Phi _{1}$ and $(2k-p,p)$ is the initial
transposition of a putative involution that is to be constructed
inductively. Terms in the second (double) summation in (\ref
{latticerecursion}) are treated similarly - they correspond to expressions $%
\gamma (\Xi )F_{N}(...)$ where $\Xi =\{m-p+1,p\}$ and $m-p+1,p$ are labels
of elements in $\Phi _{k},\Phi _{1}$ respectively (for $k\neq 1$). The third
term in (\ref{latticerecursion}) is similarly seen to coincide with $\gamma
(\Xi )F_{N}(...)$, now with $\Xi =\{p\}\subseteq $ $\Phi _{1}$. We repeat
this process in an inductive manner. It is easy to see that in this way we
construct every element of $\mathrm{Inv}(\Phi )$ exactly once, and what
emerges is the formula (\ref{Fnlattice}). This completes our discussion of
the proof of Proposition \ref{Propgennpt}. $\qed $

In order to complete our discussion of the rank one $n$-point functions we
must of course evaluate the term $F_{N}(\mathbf{1}\otimes e^{\alpha
_{1}},z_{1};\ldots ;\mathbf{1}\otimes e^{\alpha _{n}},z_{n};q)$. Before we
do that, however, it will be useful to draw some initial conclusions from
Proposition \ref{Propgennpt}. Taking all lattice elements $\alpha
_{1}=\ldots =\alpha _{n}=\beta =0$ corresponds to the case of $n$-point
functions in the rank $1$ free bosonic theory $M$. In this case all
contributions from orbits of length $1$ vanish and hence the sum over $%
\mathrm{Inv}(\Phi )$ reduces to one over $F(\Phi )$ of (\ref{FPhi}) only.
Furthermore we know that $F_{M}(\mathbf{1},z_{1};\ldots ;\mathbf{1},z_{n};q)$
is just the partition function for $M$ i.e. $1/\eta (\tau )$. Thus we arrive
at a formula for $n$-point functions for a single free boson:

\begin{corollary}
\label{gen1ptM} Let $M$ be the VOA for a single free boson. For Fock states $%
v_{1},\ldots ,v_{n}$ as in (\ref{vstate}) with corresponding labelled set $%
\Phi $ we have 
\begin{equation}
F_{M}(v_{1},z_{1};\ldots ;v_{n},z_{n};q)=\frac{1}{\eta (\tau )}\sum_{\varphi
\in F(\Phi )}\Gamma (\varphi ).  \label{BosonFM}
\end{equation}
$\qed $
\end{corollary}

The case $n=1$ of Corollary \ref{gen1ptM} was established in \cite{DMN},
where it was shown that 
\begin{equation}
Z_{M}(v;\tau )=\frac{1}{\eta (\tau )}\sum_{\varphi \in F(\Phi )}\prod
C(r,s,\tau ),  \label{FMonept}
\end{equation}
for fixed point free involutions $\varphi =\ldots (rs)\ldots $ of the
labelled set $\Phi $ labelling $v$ of (\ref{vstate}) where the product is
taken over all the transpositions $(rs)$ in $\varphi $.

An even more special, yet interesting, case arises if in Corollary \ref
{gen1ptM} we take each state $v_{i}$ to coincide with the conformal weight
one state $a$. Thus $v_{i}=a[-1].\mathbf{1}$, $1\leq i\leq n$, $\Phi $
consists of $n$ elements each carrying the label $1$, and elements of $%
\Sigma (\Phi )$ may be considered as mappings on the set $\{1,2,\ldots ,n\}$
in the usual way. If $n$ is odd then there are no fixed-point-free
involutions acting on $\Phi $, so (\ref{BosonFM}) is zero in this case. If $%
n $ is even then $\gamma (\Xi )=P_{2}(z_{ij},\tau )$ if $\Xi $ = $\Phi
_{i}\cup \Phi _{j}$ for $i\neq j$. Thus we obtain

\begin{corollary}
\label{nastatesM} Let $M$ be the VOA for a single free boson. Then for $n$
even,
\end{corollary}

\begin{equation}
F_{M}(a,z_{1};\ldots ;a,z_{n};q)=\frac{1}{\eta (\tau )}\sum_{\varphi \in
F(\Phi )}\prod P_{2}(z_{ij},\tau ),  \label{BosonFMaa}
\end{equation}
where the product ranges over the cycles of $\varphi =\ldots (ij)\ldots
\qed $

We next consider the case of an $M$-module $N$ $=M\otimes e^{\beta }$. If $%
n=1$ we necessarily have $\alpha _{1}=0$, $Z_{N}=q^{(\beta ,\beta )/2}/\eta
(\tau )$, and the labelled set $\Phi $ coincides with $\Phi _{1}$. If $%
\varphi \in \mathrm{Inv}(\Phi )$ and $\Xi =\{r\}$ is an orbit of $\varphi $
of length $1$ then from (\ref{gammaPhik}) we get $\gamma (\Xi )=\delta
_{r,1}(a,\beta )$. So $\Gamma (\varphi )$ vanishes unless all labels of the
set of fixed-points $\mathrm{Fix}(\varphi )$ are equal to $1$. In this case
we may write $\varphi =1^{|\Delta |}\varphi _{0}$ to indicate that $\varphi $
fixes a set $\Delta $ of elements with label $1$, and that $\varphi _{0}\ $%
is the fixed-point-free involution induced by $\varphi $ on the complement $%
\Phi \backslash \Delta $. We thus obtain

\begin{corollary}
\label{gen1ptN} Let $M$ be as in (\ref{nptlattice}) and let $N$ be the $M$%
-module $M\otimes e^{\beta }$. Take $v$ to be the state (\ref{vstate}) and
let $\Lambda $ denote the elements in $\Phi $ with label $1$. Then
\end{corollary}

\begin{eqnarray}
Z_{N}(v,\tau ) &=&\frac{q^{(\beta ,\beta )/2}}{\eta (\tau )}\sum_{\Delta
\subseteq \Lambda }(a,\beta )^{|\Delta |}\sum_{\varphi _{0}\in F(\Phi
\backslash \Delta )}\Gamma (\varphi _{0}),  \nonumber \\
\Gamma (\varphi _{0}) &=&\prod C(r,s,\tau ),  \label{FNoneptlattice}
\end{eqnarray}
where $\varphi _{0}=\ldots (rs)\ldots $ acting on $\Phi \backslash \Delta $. 
$\qed $

Similarly to Corollary \ref{nastatesM} we can consider again the special
case with each $v_{i}=a$, generalising (\ref{BosonFMaa}). In the above
notation we therefore have $\Lambda =\Phi $ so that

\begin{corollary}
\label{nastatesN} Let $M$ be the VOA for a single free boson with module $%
N=M\otimes e^{\beta }$. Then
\end{corollary}

\begin{equation}
F_{N}(a,z_{1};\ldots ;a,z_{n};q)=\frac{q^{(\beta ,\beta )/2}}{\eta (\tau )}%
\sum_{\Delta \subseteq \Phi }(a,\beta )^{|\Delta |}\sum_{\varphi _{0}\in
F(\Phi \backslash \Delta )}\prod P_{2}(z_{ij},\tau ),  \label{BosonFNaa}
\end{equation}
where the product ranges over the cycles of $\varphi _{0}=\ldots (ij)\ldots
\qed $

We now show how (\ref{BosonFNaa}) can be interpreted as the generator of all
free bosonic $n$-point functions for Fock states (\ref{Fockstate}). This
provides a useful insight into the structure found for these $n$-point
functions in terms of the elliptic function $P_{2}(z,\tau )$ and the scalar $%
(a,\beta )$.

\begin{proposition}
\label{Propgen} $F_{N}(a,z_{1};\ldots ;a,z_{n};q)$ is a generating function
for the $n$-point functions for all Fock states $v_{1},\ldots ,v_{n}$.
\end{proposition}

\textbf{Proof.} This follows from Lemma \ref{lemma3.1} and the expansions of 
$P_{2}$ of (\ref{P2expansion1}) and (\ref{P2expansion2}). We will illustrate
the result for $n=1$ and $n=2$. A general proof can be given along the same
lines.

From (\ref{Fnz1zn}) we obtain 
\begin{eqnarray}
F_{N}(a,z_{1};\ldots ;a,z_{n};q) &=&Z_{N}(Y[a,z_{1}]\ldots Y[a,z_{n}].%
\mathbf{1},q)  \nonumber \\
&=&\sum_{l_{1},\ldots l_{n}\in \mathbb{Z}}Z_{N}(a[-l_{1}]\ldots a[-l_{n}].%
\mathbf{1},q)z_{1}^{l_{1}-1}\ldots z_{n}^{l_{n}-1}. \nonumber\\ &&  \label{FMonepointgen}
\end{eqnarray}
The 1-point function for the bosonic Fock state $v=a[-l_{1}]\ldots a[-l_{n}].%
\mathbf{1}$ is clearly the coefficient of $z_{1}^{l_{1}-1}\ldots
z_{n}^{l_{n}-1}$ for $l_{1},\ldots ,l_{n}>0$. We then recover (\ref
{FNoneptlattice}) from the expansion for each $P_{2}(z_{ij},\tau )$ using (%
\ref{P2expansion1}).

For $n=2$ consider the $1$-point function 
\begin{equation}
Z_{N}(Y[Y[a,w_{1}]\ldots Y[a,w_{m}].\mathbf{1},w].Y[Y[a,z_{1}]\ldots
Y[a,z_{n}].\mathbf{1},z].\mathbf{1},q).  \label{FNtwopointgen}
\end{equation}
The 2-point function $F_{N}(v_{1},q_{1};v_{2},q_{2};q)$ for  $%
v_{1}=a[-l_{1}]\ldots a[-l_{m}].\mathbf{1}$ and $v_{2}=a[-k_{1}]\ldots
a[-k_{n}].\mathbf{1}$ is the coefficient of $\prod_{i=1}^{m}%
\prod_{j=1}^{n}w_{i}^{l_{i}-1}z_{j}^{m_{j}-1}$ in (\ref{FNtwopointgen}). By
associativity (\ref{Yassociativity}) and using $Y[\mathbf{1},z]=\mathrm{Id}$
eqn. (\ref{FNtwopointgen}) can be expressed as 
\[
Z_{N}(Y[a,w_{1}+w]\ldots Y[a,w_{m}+w].Y[a,z_{1}+z]\ldots Y[a,z_{n}+z].%
\mathbf{1},q). 
\]
Using (\ref{BosonFNaa}) this becomes (suppressing the $\tau$ dependence for 
clarity)
\[
\frac{q^{(\beta ,\beta )/2}}{\eta (\tau )}\sum_{\Delta \subseteq \Phi
}(a,\beta )^{|\Delta |}\sum_{\varphi _{0}\in F(\Phi \backslash \Delta
)}P_{2}(w_{ab})\ldots P_{2}(z_{cd})\ldots P_{2}(w-z+w_{e}-z_{f})\ldots 
\]
where $\varphi _{0}=(ab)\ldots (cd)\ldots (ef)\ldots \in F(\Phi \backslash
\Delta )$ with $a,b,e\ldots \in \{1,2,\ldots m\}$ and $c,d,f\ldots \in
\{1,2,\ldots n\}$. Then the coefficient of $w_{a}^{l_{a}-1}w_{b}^{l_{b}-1}$
in $P_{2}(w_{ab})$ is $C(l_{a},l_{b})$, the coefficient of $%
z_{c}^{m_{c}-1}z_{d}^{m_{d}-1}$ in $P_{2}(z_{cd})$ is $C(m_{c},m_{d})$ from (%
\ref{P2expansion1}) and the coefficient of $w_{e}^{l_{e}-1}z_{f}^{m_{f}-1}$
in $P_{2}(w-z+w_{e}-z_{f})$ is $D(l_{e},m_{f},w-z)$ from (\ref
{P2expansion2}) leading to the result (\ref{BosonFM}) in this case. $%
\qed $

We complete our discussion of bosonic $n$-point functions with two global
formulas for $1$-point functions. The first shows how to write certain $1$-point
functions with respect to $N$ in terms of $1$-point functions with respect
to $M$:

\begin{proposition}
\label{PropZNtoZM} Let notation be as in Corollary \ref{gen1ptM}. Then if $%
\varsigma $ is an indeterminate,
\end{proposition}

\begin{eqnarray}
&&Z_{N}(\exp (\sum_{m\geq 1}\frac{1}{m}a[-m]\varsigma ^{m}).\mathbf{1},\tau) \nonumber\\
&&=q^{(\beta ,\beta )/2}\exp ((a,\beta )\varsigma )Z_{M}(\exp (\sum_{m\geq 1}%
\frac{1}{m}a[-m]\varsigma ^{m}).\mathbf{1},\tau ).  \label{FNexpzeta}
\end{eqnarray}

\begin{proposition}
\label{PropZMexpPrime} Let $\lambda _{1},\ldots \lambda _{n}\,$be $n$
scalars obeying $\sum_{i=1}^{n}\lambda _{i}=0$. Then the following holds: 
\begin{equation}
Z_{M}(\exp (\sum_{m\geq 1}\frac{a[-m]}{m}\sum_{i=1}^{n}\lambda
_{i}z_{i}^{m}).\mathbf{1},\tau )=\frac{1}{\eta (\tau )}\prod_{1\leq i<j\leq
n}(\frac{K(z_{ij},\tau )}{z_{ij}})^{\lambda _{i}\lambda _{j}},
\label{FMexpzetaprime}
\end{equation}
where $K(z,\tau )$ $\,$is the prime form of (\ref{Primeform}).
\end{proposition}

We divide the proof of Proposition \ref{PropZNtoZM} into several steps.

\begin{lemma}
\label{LemZNaminusonep} Let $u\in M$ be a state such that $a[1].u=0$. For an
integer $p\geq 0$,
\end{lemma}

\begin{equation}
Z_{N}(a[-1]^{p}.u,q)=q^{(\beta ,\beta )/2}Z_{M}((a,\beta )+a[-1])^{p}.u,q).
\label{FNaminus1p}
\end{equation}

\textbf{Proof}: Suppose first that $p=0$. In this case, the lemma says that
for a state $v$ as in (\ref{vstate}) which satisfies also $e_{1}=0$, the
traces of $o(v)$ over $N$ and $M$ differ only by an overall factor of $%
q^{(\beta ,\beta )/2}$. A moment's thought shows that this follows from
Corollary \ref{gen1ptM} because the set $\Lambda $ is empty in this case.
This proves the case $p=0$ of the lemma.

We prove the general case by induction on $p$. Using lemma \ref{LemZhurec}
with $n=1$ we calculate

\begin{eqnarray*}
&&Z_{N}(a[-1]^{p+1}.u,\tau ) \nonumber\\ 
&&=(a,\beta )Z_{N}(a[-1]^{p}.u,\tau )+\sum_{k\geq
1}E_{2k}(\tau )Z_{N}(a[2k-1].a[-1]^{p}.u,\tau ) \\
&&=q^{(\beta ,\beta )/2}(a,\beta )Z_{M}((a,\beta )+a[-1])^{p}.u,\tau ) \\
&&+pE_{2}(\tau )q^{(\beta ,\beta )/2}Z_{M}(((a,\beta )+a[-1])^{p-1}.u,\tau )
\\
&&+q^{(\beta ,\beta )/2}\sum_{k\geq 2}E_{2k}(\tau )Z_{M}(a[2k-1].(a,\beta
)+a[-1])^{p}.u,\tau ),
\end{eqnarray*}
from which the result follows by using lemma \ref{LemZhurec} again. $%
\qed $

If we multiply (\ref{FNaminus1p}) by $\varsigma ^{p}$, rearrange the
constants and sum over $p$ we find

\begin{lemma}
\label{LemZNtoZM} Let notation be as above. Then
\end{lemma}

\begin{equation}
Z_{N}(\exp (a[-1]\varsigma ).u,\tau )=q^{(\beta ,\beta )/2}\exp ((a,\beta
)\varsigma )Z_{M}(\exp a[-1]\varsigma ).u,\tau ).
\end{equation}
$\qed $

Choosing $u=\exp (\sum_{m\geq 2}\frac{a[-m]}{m}\varsigma ^{m}).\mathbf{1\,}$%
we note that Proposition \ref{PropZNtoZM} is a special case of lemma \ref
{LemZNtoZM}.

\textbf{Proof of Proposition \ref{PropZMexpPrime}}: When expanded as a sum,
the left-hand-side of (\ref{FMexpzetaprime}) can be written in the form

\begin{equation}
\sum_{v}\prod_{k=1}^{p}\frac{1}{e_{k}!}(\frac{1}{k}\sum_{i=1}^{n}\lambda
_{i}z_{i}^{k})^{e_{k}}Z_{M}(v,\tau ),  \label{FMexpb}
\end{equation}
where $v$ ranges over the basis elements (\ref{vstate}). In this regard one
should note that as a consequence of Corollary \ref{gen1ptM}, $Z_{M}(v,\tau
)=0$ whenever $\sum_{k=1}^{p}ke_{k}$ is odd. The argument that we use to
establish the equality (\ref{FMexpzetaprime}) involves the use of a
combinatorial technique that proliferates in other work on genus two and
higher VOAs \cite{MT}.

Use the case $N=M$ of (\ref{FNoneptlattice}) and use Corollary \ref{gen1ptM}
to write (\ref{FMexpb}) as

\begin{equation}
\frac{1}{\eta (\tau )}\sum_{v}\prod_{k=1}^{p}\frac{1}{e_{k}!}\sum_{\varphi
\in F(\Phi _{v})}\prod_{k=1}^{p}\Gamma (\varphi )(\frac{1}{k}%
\sum_{i=1}^{n}\lambda _{i}z_{i}^{k})^{e_{k}}  \label{efactgamma}
\end{equation}
where $\Phi _{v}$ is the labelled set determined by $v$. Fix for a moment an
element $\varphi $, considered as a product of transpositions. We can
represent $\varphi $ by a graph with nodes labelled by positive integers
corresponding to the elements of $\Phi _{v}$, two nodes being connected
precisely when $\varphi $ interchanges the nodes in question. Such a graph
masquerades under various names in combinatorics: bipartite graph, or a 
\textit{complete matching} (cf. \cite{LW}, for example), and is nothing more
than another way to think about fixed-point-free involutions. Pictorially
the complete matching looks like

\begin{center}
\[
\begin{array}{c}
r_{1}\bullet \text{-----}\bullet s_{1} \\ 
r_{2}\bullet \text{-----}\bullet s_{2} \\ 
\vdots \\ 
r_{b}\bullet \text{-----}\bullet s_{b} \\ 
\mathrm{Fig.1}
\end{array}
\]
\end{center}

\smallskip We denote the complete matching determined by $\varphi $ by the
symbol $\mu _{\varphi }$. Any complete matching on a set labelled by
positive integers corresponds to a fixed-point-free involution and a state $%
v $. Let us agree that two complete matchings $\mu _{1},\mu _{2}$ are 
\textit{isomorphic} if there is a bijection from the node set of $\mu _{1}$
to the node set of $\mu _{2}$ that preserves labels. We may identify the
node sets of the two matchings, call it $\Phi $, in which case an
isomorphism may be realized by an element in the symmetric group $\Sigma
(\Phi )$. More precisely, let us define the \textit{label subgroup} of $%
\Sigma (\Phi )$ to be the subgroup $\Lambda (\Phi )$ consisting of all
permutations of $\Phi $ that preserve labels. Then it is the case that an
isomorphism between $\mu _{1}$ and $\mu _{2}$ may be realized by an element
in the label subgroup. Note that our notation implies that $|\Lambda (\Phi
)|=\prod e_{k}!$.

Now consider (\ref{efactgamma}). The expression to the right of the second
summation, which we denote by 
\begin{equation}
\hat{\Gamma}(\mu _{\varphi })=\prod_{k=1}^{p}\Gamma (\varphi )(\frac{1}{k}%
\sum_{i=1}^{n}\lambda _{i}z_{i}^{k})^{e_{k}}  \label{Gammahat}
\end{equation}
is determined by the complete matching $\mu _{\varphi }$. By what we have
said, every complete matching on a set labelled by positive integers occurs
in this way as $v$ ranges over the preferred Fock basis of $M$ (on the
cylinder) and $\varphi $ ranges over $F(\Phi _{v})$, while the factor $\prod 
\frac{1}{e_{k}!}$ may be interpreted as averaging over all complete
matchings with a given set of labels. The only duplication that occurs is
due to \textit{automorphisms} of the complete matching. The upshot is that
expression (\ref{efactgamma}) is equal to 
\begin{equation}
\frac{1}{\eta (\tau )}\sum_{\mu }\frac{\hat{\Gamma}(\mu )}{|\mathrm{Aut}(\mu
)|}  \label{autinvgammahat}
\end{equation}
where $\mu $ ranges over all \textit{isomorphism classes} of complete
matchings labelled by positive integers.

For a given complete matching $\mu $ as in Fig. 1, let $E=E(r,s)$ denote an
edge with labels $r,s$, and let $m(E)$ denote the \textit{multiplicity} of $%
E $ in $\mu $. Thus we may represent $\mu $ symbolically by its
decomposition $\mu =\sum m(E)E$ into isomorphism classes of labelled edges.
Now it is evident that there is an isomorphism of groups

\[
\mathrm{Aut}(\mu )\simeq \prod_{E}\mathrm{Aut}(E)\wr \Sigma _{m(E)}, 
\]
a direct product, indexed by isomorphism classes of labelled edges, of
groups which are themselves the (regular) wreathed product of $\mathrm{Aut}%
(E)$ and $\Sigma _{m(E)}$. In particular, we have 
\begin{equation}
|\mathrm{Aut}(\mu )|=\prod_{E}m(E)!|\mathrm{Aut}(E)|^{m(E)}.
\label{orderautmu}
\end{equation}
Note that $|\mathrm{Aut}(E)|\leq 2$, with equality only if the two node
labels of $E$ are equal. Next it is easy to see that the expression $\hat{%
\Gamma}(\mu )$ is \textit{multiplicative} \textit{over edges}. In other
words, we have 
\begin{equation}
\hat{\Gamma}(\mu )=\prod_{E}\hat{\Gamma}(E)^{m(E)},  \label{gammahatmu}
\end{equation}
and for an edge $E(r,s)$ we have 
\begin{eqnarray}
\hat{\Gamma}(E) &=&C(r,s,\tau )(\frac{1}{r}\sum_{i=1}^{n}\lambda
_{i}z_{i}^{r})(\frac{1}{s}\sum_{j=1}^{n}\lambda _{j}z_{j}^{s})  \nonumber \\
&=&\frac{(-1)^{r+1}}{r+s}\binom{r+s}{s}E_{r+s}(\tau )\sum_{i=1}^{n}\lambda
_{i}z_{i}^{r}\sum_{j=1}^{n}\lambda _{j}z_{j}^{s}.  \label{gammahatE}
\end{eqnarray}
Substitute (\ref{orderautmu}) and (\ref{gammahatmu}) in (\ref{autinvgammahat}%
) to obtain the expression 
\begin{equation}
\frac{1}{\eta (\tau )}\prod_{E}\exp (\frac{\hat{\Gamma}(E)}{|\mathrm{Aut}(E)|%
})=\frac{1}{\eta (\tau )}\exp (\sum_{E_{or}}\frac{\hat{\Gamma}(E_{or})}{2}),
\label{expgammahatE}
\end{equation}
where $E$ ranges over all labelled edges $r \bullet$------$\bullet s$
and $E_{or}$ ranges over all \textit{oriented} edges 
$r\bullet $--->---$\bullet s$ (which have trivial automorphism group). Using (%
\ref{gammahatE}), the expression (\ref{expgammahatE}) is in turn equal to

\begin{equation}
\frac{1}{\eta (\tau )}\exp (\frac{1}{2}\sum_{k\geq 1}E_{2k}(\tau )\frac{1}{2k%
}\sum_{r=0}^{2k}(-1)^{r+1}\binom{2k}{r}\sum_{i=1}^{n}\lambda
_{i}z_{i}^{r}\sum_{j=1}^{n}\lambda _{j}z_{j}^{2k-r}).  \label{exptwolambda}
\end{equation}
But, 
\[
\frac{1}{2}\sum_{r=0}^{2k}(-1)^{r+1}\binom{2k}{r}\sum_{i=1}^{n}\lambda
_{i}z_{i}^{r}\sum_{j=1}^{n}\lambda _{j}z_{j}^{2k-r}=-\sum_{1\leq i<j\leq
n}\lambda _{i}\lambda _{j}z_{ij}^{2k}. 
\]
and hence using (\ref{P0}), (\ref{Primeform}) we find (\ref{exptwolambda})
becomes 
\[
\frac{1}{\eta (\tau )}\prod_{1\leq i<j\leq n}\exp (-\lambda _{i}\lambda
_{j}(P_{0}(z_{ij},\tau )+\log z_{ij}))=\frac{1}{\eta (\tau )}\prod_{1\leq
i<j\leq n}[\frac{K(z_{ij},\tau )}{z_{ij}}]^{\lambda _{i}\lambda _{j}}. 
\]
This completes the proof of Proposition \ref{PropZMexpPrime}. $\qed $

We now turn our attention on the second factor $F_{N}$ of (\ref{QNFN}) where
we abbreviate $\mathbf{1}\otimes e^{\alpha _{i}}\,$ by $e^{\alpha _{i}}$
again.

\begin{proposition}
\label{Propnptlattice} Let $M$ and $N$ be as above, and let $\alpha
_{1},\ldots ,\alpha _{n}$ be lattice elements in the rank one even lattice $%
L $ satisfying $\alpha _{1}+\ldots +\alpha _{n}=0$. Then
\end{proposition}

\begin{eqnarray}
&&F_{N}(e^{\alpha _{1}},z_{1};...;e^{\alpha _{n}},z_{n};q)  \nonumber \\
&=&\frac{q^{(\beta ,\beta )/2}}{\eta (\tau )}\prod_{1\leq r\leq n}\exp
((\beta ,\alpha _{r})z_{r})\prod_{1\leq i<j\leq n}\epsilon (\alpha
_{i},\alpha _{j})K(z_{ij},\tau )^{(\alpha _{i},\alpha _{j})}.
\label{FNexpalpha}
\end{eqnarray}

\textbf{Proof}: Use (\ref{Fnz1zn}) of Lemma \ref{lemma3.1} to rewrite the
LHS of (\ref{FNexpalpha}) as

\begin{equation}
Z_{N}(o(Y[e^{\alpha _{1}},z_{1}]\ldots Y[e^{\alpha _{n}},z_{n}].\mathbf{1}%
);q)  \label{ZNalpha}
\end{equation}
Referring to (\ref{Yealpha}), by repeated use of the identity 
\[
Y_{+}(e^{\alpha },z)Y_{-}(e^{\beta },w)=(\frac{z-w}{z})^{(\alpha ,\beta
)}Y_{-}(e^{\beta },w)Y_{+}(e^{\alpha },z) 
\]
(for $|z|>|w|$\thinspace ) and using (\ref{ealpha}) to (\ref{cocycleproduct}%
) we find (\ref{ZNalpha}) is 
\begin{equation}
\prod_{1\leq r<s\leq n}z_{rs}{}^{(\alpha _{r},\alpha _{s})}\epsilon (\alpha
_{r},\alpha _{s})Z_{N}(\exp (\sum_{m>0}\sum_{i=1}^{n}\frac{\alpha _{i}[-m]}{m%
}z_{i}^{m}).\mathbf{1,}\tau ).  \label{ZNalpha2}
\end{equation}
The operator corresponding to $m=1$ in the exponential in (\ref{ZNalpha2})
may be written in the form $a[-1]\varsigma $ where $\varsigma =$ $%
\sum_{k=1}^{n}(a,\alpha _{k})z_{k}$ so that, from Lemma \ref{LemZNtoZM},

\begin{eqnarray*}
&&Z_{N}(\exp (\sum_{m>0}\sum_{i=1}^{n}\frac{\alpha _{i}[-m]}{m}z_{i}^{m}).%
\mathbf{1,}\tau ) \\
&=&q^{(\beta ,\beta )/2}\prod_{i=1}^{n}\exp ((\beta ,\alpha
_{i})z_{i})Z_{M}(\exp (\sum_{m>0}\sum_{i=1}^{n}\frac{\alpha _{i}[-m]}{m}%
z_{i}^{m}).\mathbf{1,}\tau )
\end{eqnarray*}
Now use Proposition \ref{PropZMexpPrime} with $\lambda _{i}=(a,\alpha _{i})$
so that we find 
\[
Z_{M}(\exp (\sum_{m>0}\sum_{i=1}^{n}\frac{\alpha _{i}[-m]}{m}z_{i}^{m}).%
\mathbf{1,}\tau )=\frac{1}{\eta (\tau )}\prod_{1\leq i<j\leq n}[\frac{%
K(z_{ij},\tau )}{z_{ij}}]^{(\alpha _{i},\alpha _{j})}. 
\]
This completes the proof of the Proposition. $\qed $

\smallskip We now consider the lattice VOA $V_{L}$ constructed from a rank $%
l $ even lattice $L$ as described in section 2. We recall that $%
a_{1},a_{2},...a_{l}$ is an orthonormal basis for $\mathfrak{H}$ with respect to
the non-degenerate symmetric bilinear form $(,)$. We let $M$ be the rank $l$
Heisenberg vertex operator algebra and let $N=M\otimes e^{\beta }$ be a
simple $M$-module with $\beta \in L$, $h=(\beta ,\beta )/2$ the conformal
weight of the highest weight vector of $N$. Then $M\simeq M^{1}\otimes
\ldots \otimes M^{l}$, $\,$the tensor product of $l$ copies of the rank $1$
Heisenberg VOA, and is spanned by the Fock states

\begin{equation}
v=a_{1}[-1]^{e_{1}}\ldots a_{1}[-p]^{e_{p}}\ldots a_{l}[-1]^{f_{1}}\ldots
a_{l}[-q]^{f_{q}}.\mathbf{1}  \label{genvstate}
\end{equation}
where $e_{1},\ldots ,f_{q}\,$are non-negative integers.

We now give a general closed formula for all rank $l$ lattice $n$-point
functions (\ref{nptlattice}) where $v_{1},\ldots ,v_{n}$ are Fock states of
the form (\ref{genvstate}). Viewing each vector $v_{i}$ as an element of $%
M^{1}\otimes \ldots \otimes M^{l}$ we define $\Phi _{i}^{r}\,\,\,$as the
labelled set for the $r^{\mathrm{th}}$ tensored vector of $v_{i}$ e.g. $\Phi
_{i}^{1}$ contains $\,1$ with multiplicity $e_{1}$, $2$ with multiplicity $%
e_{2}$ etc. Therefore, $v_{i}$ is determined by the labelled set $\Phi
_{i}=\Phi _{i}^{1}\cup \ldots \cup \Phi _{i}^{l}$, the disjoint union of $l$
sets. We also define $\Phi ^{r}=\bigcup_{1\leq i\leq n}\Phi _{i}^{r}$ to be
the labelled set for the $r^{\mathrm{th}}$ tensored vectors of the $n$
vectors $v_{1},\ldots ,v_{n}$. Then we have:

\begin{theorem}
\label{Big Theorem} Let $v_{1},\ldots ,v_{n}$ be states of the form (\ref
{genvstate}) in the rank $l$ free boson theory $M$ and let $\Phi ^{1},\ldots
\Phi ^{l}$ be the labelled sets defined as above by these states. Then the
following holds for lattice elements $\alpha _{1},\ldots ,\alpha _{n}\in L$
satisfying $\alpha _{1}+\ldots +\alpha _{n}$ $=0$: 
\begin{eqnarray}
&&F_{N}(v_{1}\otimes e^{\alpha _{1}},z_{1};\ldots ;v_{n}\otimes e^{\alpha
_{n}},z_{n};q)  \nonumber \\
&=&Q_{N}(v_{1},z_{1};\ldots ;v_{n},z_{n};q)F_{N}(e^{\alpha _{1}},z_{1};\ldots
;e^{\alpha _{n}},z_{n};q),  \label{FnBigtheorem}
\end{eqnarray}
where 
\begin{equation}
Q_{N}(v_{1},z_{1};\ldots ;v_{n},z_{n};q)=\prod_{1\leq r\leq l}\sum_{\varphi
^{r}\in \mathrm{Inv}(\Phi ^{r})}\Gamma (\varphi ^{r}),  \label{GenQn}
\end{equation}
and 
\begin{eqnarray}
&&F_{N}(e^{\alpha _{1}},z_{1};\ldots ;e^{\alpha _{n}},z_{n};q)  \nonumber \\
&=&\frac{q^{(\beta ,\beta )/2}}{\eta (\tau )^{l}}\prod_{1\leq r\leq n}\exp
((\beta ,\alpha _{r})z_{r})\prod_{1\leq i<j\leq n}\epsilon (\alpha
_{i},\alpha _{j})K(z_{ij},\tau )^{(\alpha _{i},\alpha _{j})}.
\label{Genlattice}
\end{eqnarray}
\end{theorem}

\textbf{Proof.} We sketch the proof which follows along the same lines as
the rank one case described above. We firstly apply Lemma \ref{LemZhurec} to
the $r^{\mathrm{th}}$ tensored vectors labelled by $\Phi ^{r}$. Following
the same argument for the rank one case in Proposition \ref{Propgennpt},
this results in (\ref{FnBigtheorem}) and (\ref{GenQn}). We secondly evaluate
the LHS of (\ref{Genlattice}) as in Proposition \ref{Propnptlattice} using
Proposition \ref{PropZMexpPrime} with $\lambda _{i}^{r}=(a_{r},\alpha _{i})$
for $r=1,\ldots ,l$ to obtain (\ref{Genlattice}). $\qed $

We conclude this section with the first non-trivial examples of lattice $n$%
-point functions which occur for $n=2$. From Theorem \ref{Big Theorem} and
recalling (\ref{cocycleunity}) we have:

\begin{corollary}
\label{Cor3.13} For a rank $l$ lattice theory with $N$ as above and with
states $e^{\alpha }$ and $e^{-\alpha }$ we have:
\end{corollary}

\begin{equation}
F_{N}(e^{\alpha },z_{1};e^{-\alpha },z_{2};q)=\frac{q^{(\beta ,\beta )/2}}{%
\eta ^{l}(\tau )}\frac{\exp ((\beta ,\alpha )z_{12})}{K(z_{12},\tau
)^{(\alpha ,\alpha )}}.  \label{lattice2pt}
\end{equation}
$\qed $

Taking the sum over all $\beta \in L$ we immediately obtain:

\begin{corollary}
For $V=V_{L}$, the lattice vertex operator algebra for a rank $l$ even
lattice $L$, and for states $e^{\alpha }$ and $e^{-\alpha }$ we have:

\begin{equation}
F_{V_{L}}(e^{\alpha },z_{1};e^{-\alpha },z_{2};q)=\frac{1}{\eta ^{l}(\tau )}%
\frac{\Theta _{\alpha ,L}(\tau ,z_{12}/2\pi i)}{K(z_{12},\tau )^{(\alpha
,\alpha )}},
\end{equation}
where 
\begin{equation}
\Theta _{\alpha ,L}(\tau ,z)=\sum_{\beta \in L}\exp (2\pi i[\frac{(\beta
,\beta )}{2}\tau +(\beta ,\alpha )z]),
\end{equation}
is a Jacobi form of weight $l\,$ and index $(\alpha ,\alpha )/2$ \cite{EZ}. $%
\qed $
\end{corollary}

\section{The Elliptic Properties of $n$-point Functions}

In this section we will consider the elliptic properties of the $n$-point
functions described in the previous section. For vertex operator algebras
satisfying the so-called $C_{2}$ condition, Zhu has shown that every $n$%
-point function is meromorphic and periodic in each parameter $z_{i}\,$and
is therefore elliptic \cite{Z}. The $C_{2}$ condition does not hold for
simple modules of free bosonic theories but nevertheless all $n$-point
functions (\ref{FnBigtheorem}) are found to be meromorphic and either
elliptic for the free bosonic $n$-point functions or quasi-elliptic for the
lattice $n$-point functions. In this section we will consider these $n$%
-point functions from first principles where our aim is to provide further
insight into the structure found for these functions. In particular, we will
show how (\ref{BosonFNaa}), the generating function for free bosonic $n$%
-point functions and the lattice $n$-point function (\ref{FNexpalpha}) are
the unique elliptic or quasi-elliptic functions determined by permutation
symmetry, periodicity and certain natural singularity and normalisation
properties.

We begin with a general statement about $n$-point functions:

\begin{lemma}
\label{Lemma4.1} The $n$-point function $F_{N}=F_{N}(v_{1}\otimes e^{\alpha
_{1}},z_{1};\ldots ;v_{n}\otimes e^{\alpha _{n}},z_{n};q)$ for $\alpha
_{1}+\ldots +\alpha _{n}$ $=0$ obeys the following:

\begin{description}
\item[(i)]  $F_{N}$ is symmetric under all permutations of its indices.

\item[(ii)]  $F_{N}$ is a function of $z_{ij}=$ $z_{i}-z_{j}$.

\item[(iii)]  $F_{N}$ is non-singular at $z_{ij}$ $\neq 0$ for all $i\neq j$.

\item[(iv)]  $F_{N}$ is periodic in $z_{i}$ with period $2\pi i$.

\item[(v)]  $F_{N}$ is quasi-periodic in $z_{i}$ with period $2\pi i\tau $
and multiplier 
\begin{equation}
q^{(\alpha _{i},\alpha _{i})/2+(\alpha _{i},\beta )}q_{i}^{(\alpha
_{i},\alpha _{i})}.  \label{quasiperiodmult}
\end{equation}
\end{description}
\end{lemma}

\textbf{Proof}. (i) Apply the general locality property for vertex operators
e.g. \cite{Ka} 
\[
(z-w)^{N}Y(u,z).Y(v,w)=(z-w)^{N}Y(v,w).Y(u,z), 
\]
for $N$ sufficiently large, to all adjacent pairs of operators in (\ref
{Fnz1zn}) of Lemma \ref{lemma3.1}.

(ii) This follows from (i) and (\ref{Fnziminuszn}) of Lemma \ref{lemma3.1}.

(iii) Suppose that $F_{N}$ is singular at $z_{n}=z_{0}$ for some $z_{0}\neq $
$z_{j}$ for all $j=1,\ldots ,n-1$. Using (ii) we may assume that $z_{0}=0$
by redefining $z_{i}$ to be $z_{i}-z_{0}$ for all $i$. But $F_{N}$ cannot be
singular at $z_{n}=0$ from (\ref{Fnz1zn}) of Lemma \ref{lemma3.1} since $%
Y[v_{n}\otimes e^{\alpha _{n}},z_{n}].\mathbf{1}|_{z_{n}=0}=v_{n}$ and hence
the result follows.

(iv) This follows from the integrality of conformal weights.

(v) Using (i) we have 
\begin{eqnarray}
F_{N}&=&q^{-c/24}Tr_{N}Y(q_{2}^{L(0)}v_{2}\otimes e^{\alpha _{2}},q_{2}).\nonumber\\
&&\ldots
Y(q_{n}^{L(0)}v_{n}\otimes e^{\alpha _{n}},q_{n}).Y(q_{1}^{L(0)}v_{1}\otimes
e^{\alpha _{1}},q_{1})q^{L(0)}.\nonumber 
\end{eqnarray}
Under $z_{1}\rightarrow z_{1}+2\pi i\tau \,$ $\,$we have $F_{N}\rightarrow 
\hat{F}_{N}$ where 
\begin{eqnarray}
\hat{F}_{N}&=&q^{-c/24}Tr_{N}Y(q_{2}^{L(0)}v_{2}\otimes e^{\alpha
_{2}},q_{2}).\nonumber\\
&&\ldots Y(q_{n}^{L(0)}v_{n}\otimes e^{\alpha
_{n}},q_{n}).q^{L(0)}.Y(q_{1}^{L(0)}v_{1}\otimes e^{\alpha _{1}},q_{1}),
\label{Fnhat}
\end{eqnarray}
using (\ref{L0scaling}). Consider the co-cycle parts of the vertex operators
within $\hat{F}_{N}$. Using (\ref{ealphacomm}) to (\ref{zalpha}) we see that 
\begin{eqnarray}
&&\prod_{2\leq r\leq n}e^{\alpha _{r}}.q_{r}^{\alpha
_{r}}.q^{L(0)}.e^{\alpha _{1}}.q_{1}^{\alpha _{1}}.(v\otimes e^{\beta }) 
\nonumber \\
&=&q^{(\alpha _{1},\alpha _{1})/2+(\alpha _{1},\beta )}q_{1}^{(\alpha
_{1},\alpha _{1})}\prod_{1\leq r\leq n}e^{\alpha _{r}}.q_{r}^{\alpha
_{r}}.q^{L(0)}.(v\otimes e^{\beta }).  \label{cocycleL0}
\end{eqnarray}
Since $N=\bigoplus_{n\in \mathbb{Z}}M_{n}\otimes e^{\beta }$, the graded trace (%
\ref{Fnhat}) decomposes into finite dimensional traces over $M_{n}$ to which
we may apply the standard trace property $TrAB=TrBA$ on the remaining parts
of the vertex operators within $\hat{F}_{N}$. Hence we obtain 
\[
\hat{F}_{N}=q^{(\alpha _{1},\alpha _{1})/2+(\alpha _{1},\beta
)}q_{1}^{(\alpha _{1},\alpha _{1})}F_{N}, 
\]
as required. By (i) we obtain the quasi-periodicity (\ref{quasiperiodmult})
for each $z_{i}$. $\qed $

Let us now consider the generating function, $F_{N}(a,z_{1};\ldots
;a,z_{n};q)$, for all Fock state $n$-point functions for the rank one case.
Note that $F_{N}(q)\equiv Z_{N}(q)=q^{(\beta ,\beta )/2}/\eta (\tau )$ for $%
n=0$. From Lemma \ref{Lemma4.1}, $F_{N}(a,z_{1};\ldots ;a,z_{n};q)$ $\,$is
periodic in each $z_{i}$ with periods $2\pi i\,$and $2\pi i\tau $. The
singularity structure at $z_{ij}=0$ is determined by

\begin{lemma}
\label{Lemma4.3} For $n\geq 2$ and for $i\neq j$, $F_{N}(a,z_{1};\ldots
;a,z_{n};q)$ has the following leading behaviour in its (formal) Laurent
expansion in $z_{ij}$ 
\begin{equation}
F_{N}(a,z_{1};\ldots ;a,z_{n};q)=\frac{1}{z_{ij}^{2}}F_{N}(a,z_{1};\ldots ;%
\hat{a},\hat{z}_{i}\ldots ;\hat{a},\hat{z}_{j}\ldots ;a,z_{n};q)+\ldots ,
\label{FNa1tonres}
\end{equation}
where $\hat{a},\hat{z}_{i}$ and $\hat{a},\hat{z}_{j}$ denotes the deletion
of the corresponding vertex operators resulting in an $n-2$ point function.
\end{lemma}

\textbf{Proof.} Using Lemma \ref{Lemma4.1} (i) it suffices to consider the
expansion in $z_{n-1n}$. The result then follows from (\ref{Fnziminuszn}) of
Lemma \ref{lemma3.1} and using 
\[
Y[a,z_{n-1n}].a=\frac{1}{z_{n-1n}^{2}}\mathbf{1}+\sum_{k\geq
1}z_{n-1n}^{k-1}a[-k].a. 
\]

\smallskip $\qed $

We also have the following integral normalisation:

\begin{lemma}
\label{Lemma4.4} For each $i=1,\ldots ,n\geq 1$ 
\begin{equation}
\frac{1}{2\pi i}\int_{0}^{2\pi i}F_{N}(a,z_{1};\ldots
;a,z_{n};q)dz_{i}=(a,\beta )F_{N}(a,z_{1};\ldots ;\hat{a},\hat{z}_{i};\ldots
;a,z_{n};q),  \label{aPeriod}
\end{equation}
where $\hat{a},\hat{z}_{i}$ denotes the deletion of the corresponding vertex
operator giving an $n-1$ point function.
\end{lemma}

\textbf{Proof.} Using Lemma \ref{Lemma4.1} (i) it suffices to consider $i=1$
only. Then the integral is 
\begin{eqnarray*}
&&\frac{1}{2\pi i}\int_{0}^{2\pi i}F_{N}(a,z_{1};\ldots ;a,z_{n};q)dz_{1} \\
&=&Tr_{N}\frac{1}{2\pi i}\oint_{\mathcal{C}%
_{1}}Y(a,q_{1})dq_{1}Y(q_{2}a,q_{2})\ldots Y(q_{n}v_{n},q_{n})q^{L(0)-1/24}
\\
&=&Tr_{N}o(a)Y(q_{2}a,q_{2})\ldots Y(q_{n}v_{n},q_{n})q^{L(0)-1/24} \\
&=&(a,\beta )Tr_{N}Y(q_{2}a,q_{2})\ldots Y(q_{n}v_{n},q_{n})q^{L(0)-1/24}.
\end{eqnarray*}

where $\mathcal{C}_{1}\,$denotes a closed contour surrounding $q_{1}=0$. $%
\qed $

We now show that $F_{N}(a,z_{1};\ldots ;a,z_{n};q)$ $\,$is uniquely
determined to be given by (\ref{BosonFNaa}) of Corollary \ref{nastatesN} as
follows:

\begin{proposition}
\label{Prop4.5} $F_{N}(a,z_{1};\ldots ;a,z_{n};q)\,$ is the unique
meromorphic function in $z_{i}\in \mathbb{C}/\{2\pi i(m+n\tau )|m,n\in \mathbb{Z}%
\} $ obeying Lemmas \ref{Lemma4.1}, \ref{Lemma4.3} and \ref{Lemma4.4} and is
given by (\ref{BosonFNaa}).
\end{proposition}

\textbf{Proof.} If $F_{N}(a,z_{1};\ldots ;a,z_{n};q)\,$ is meromorphic then
it is an elliptic function from Lemma \ref{Lemma4.1} (iv) and (v). We prove
the required result by induction. For $n=1$, $F_{N}(a,z_{1};q)$ has no poles
from Lemma \ref{Lemma4.1} (iii) and is therefore constant in $z_{1}$. Then (%
\ref{aPeriod}) of Lemma \ref{Lemma4.4} implies 
\[
F_{N}(a,z_{1};q)=(a,\beta )\frac{q^{(\beta ,\beta )/2}}{\eta (\tau )}, 
\]
in agreement with the RHS of (\ref{BosonFNaa}) in this case.

Next consider the elliptic function 
\begin{eqnarray}
G(z_{1},\ldots ,z_{n})&=&F_{N}(a,z_{1};\ldots
;a,z_{n};q)\nonumber\\
&&-\sum_{i=2}^{n}P_{2}(z_{1i})F_{N}(a,z_{2};\ldots ;\hat{a},\hat{z}%
_{i};\ldots ;a,z_{n};q), \nonumber
\end{eqnarray}
where $\hat{a},\hat{z}_{i}$ denotes the deletion of the corresponding vertex
operator. From Lemma \ref{Lemma4.3}, $G$ is holomorphic and elliptic in $%
z_{1}$ and is therefore independent of $z_{1}$. Next note that although not
an elliptic function, $P_{1}(z,\tau )$ is periodic with period $2\pi i$ and
so $\int_{0}^{2\pi i}P_{2}(z_{1i},\tau )dz_{1}=0$. Hence Lemma \ref{Lemma4.4}
implies 
\begin{eqnarray}
F_{N}(a,z_{1};\ldots ;a,z_{n};q) &=&(a,\beta )F_{N}(a,z_{2};\ldots
;a,z_{n};q)  \nonumber \\
&&+\sum_{i=2}^{n}P_{2}(z_{1i})F_{N}(a,z_{2};\ldots ;\hat{a},\hat{z}%
_{i};\ldots ;a,z_{n};q). \nonumber\\
&& \label{Fnrecurrence}
\end{eqnarray}
But this recurrence relation is precisely (\ref{FNa1Zhu}) with $v_{1}=%
\mathbf{1},v_{2}=\ldots =v_{n}=a$ and $\alpha _{i}=0$. Thus we obtain the
RHS of (\ref{BosonFNaa}) by induction as before. $\qed $

\begin{remark}
The rank $l$ result follows as before by considering the tensor product of $%
l $ rank one Heisenberg VOAs.
\end{remark}

We now consider the lattice $n$-point function $F_{N}(e^{\alpha
_{1}},z_{1};...;e^{\alpha _{n}},z_{n};q)$ $\,$for a rank $l$ lattice. We
firstly note the following:

\begin{lemma}
\label{Lemmaalpharec} For $n\geq 1$ and for $i\neq j$, the (formal) Laurent
expansion in $z_{ij}$ of $F_{N}(e^{\alpha _{1}},z_{1};...;e^{\alpha
_{n}},z_{n};q)$ has leading behaviour 
\begin{eqnarray}
&&F_{N}(e^{\alpha _{1}},z_{1};...;e^{\alpha _{n}},z_{n};q)  \nonumber \\
&&=\epsilon (\alpha _{i},\alpha _{j})z_{ij}^{(\alpha _{i},\alpha
_{j})}F_{N}(e^{\alpha _{1}},z_{1};...\hat{e}^{\alpha _{i}},\hat{z}%
_{i};\ldots e^{\alpha _{i}+\alpha _{j}},z_{j};\ldots ;e^{\alpha
_{n}},z_{n};q)+\ldots ,\nonumber \\
&&
\end{eqnarray}
where $\hat{e}^{\alpha _{i}},\hat{z}_{i}$ denotes the deletion of the
corresponding vertex operator resulting in an $n-1$ point function.
\end{lemma}

\textbf{Proof.} Using Lemma \ref{Lemma4.1} (i) it suffices to consider the
expansion in $z_{n-1n}$. The result then follows from (\ref{ealpha}), (\ref
{zalpha}) and (\ref{Fnziminuszn}) of Lemma \ref{lemma3.1} to find 
\[
Y[e^{\alpha _{n-1}},z_{n-1n}].e^{\alpha _{n}}=\epsilon (\alpha _{n-1},\alpha
_{n})z_{n-1n}^{(\alpha _{n-1},\alpha _{n})}e^{\alpha _{n}+\alpha
_{n-1}}+\ldots 
\]

\smallskip $\qed $

Next recall the following properties for the prime form $K(z,\tau )$ e.g. 
\cite{Mu}

\begin{lemma}
\label{LemmaPrime} The genus one prime form $K(z,\tau )$ is a holomorphic
function on $\mathbb{C}/\{2\pi i(m+n\tau )|m,n\in \mathbb{Z}\}$ given by 
\begin{eqnarray}
K(z,\tau )&=&-\frac{i\theta _{1}(z,\tau )}{\eta (\tau )^{3}},
\label{Kthetaeta} \\
\theta _{1}(z,\tau )&\equiv& \sum_{n\in \mathbb{Z}}\exp (\pi i\tau
(n+1/2)^{2}+(n+1/2)(z+i\pi )),  \label{theta1}
\end{eqnarray}
where $K(z,\tau )$ is quasi-periodic in $z$ with period $2\pi i\,$ and
multiplier $-1$ and with period $2\pi i\tau $ and multiplier $%
-q^{-1/2}q_{z}^{-1}$. Furthermore $K(z,\tau )$ has a unique zero at $z=0$ on 
$\mathbb{C}/\{2\pi i(m+n\tau )|m,n\in \mathbb{Z}\}$. $\qed $
\end{lemma}

We finally show that $F_{N}(e^{\alpha _{1}},z_{1};...;e^{\alpha
_{n}},z_{n};q)$ is uniquely determined to be given by (\ref{Genlattice}) as
follows:

\begin{proposition}
\label{Propalpha} $F_{N}(e^{\alpha _{1}},z_{1};...;e^{\alpha _{n}},z_{n};q)\,
$ is the unique meromorphic function in $z_{i}\in \mathbb{C}/\{2\pi i(m+n\tau
)|m,n\in \mathbb{Z}\}$ obeying Lemmas \ref{Lemma4.1} and \ref{Lemmaalpharec}.
\end{proposition}

\textbf{Proof.} If $F_{N}(e^{\alpha _{1}},z_{1};...;e^{\alpha
_{n}},z_{n};q)\,$ is meromorphic then consider the meromorphic function 
\begin{equation}
G(z_{1},\ldots ,z_{n})=\frac{F_{N}(e^{\alpha _{1}},z_{1};...;e^{\alpha
_{n}},z_{n};q)}{\prod_{1\leq r\leq n}\exp ((\beta ,\alpha
_{r})z_{r})\prod_{1\leq i<j\leq n}\epsilon (\alpha _{i},\alpha
_{j})K(z_{ij},\tau )^{(\alpha _{i},\alpha _{j})}}.  \label{Galpha}
\end{equation}
We wish to show that $G=F_{N}(q)=\eta (\tau )^{l}$. We prove this by
induction in $n$. It is easy to see that $G$ is periodic with periods $2\pi
i $, $2\pi i\tau $ using Lemma \ref{Lemma4.1} (iv),(v) and Lemma \ref
{LemmaPrime} and hence $G$ is elliptic in $z_{i}.$ $G$ is also a function of 
$z_{ij}$ and considering the Laurent expansion in $z_{ij}$ one finds that
the leading term is given by

\[
G(z_{1},\ldots ,z_{n})=F_{N}(q)+\ldots 
\]
using Lemma \ref{Lemmaalpharec} and (\ref{cocycleproduct}) and induction.
Hence $G$ is regular at $z_{ij}=0$. But the denominator of $G$ has zeros
possible only at $z_{ij}=0$ (for $(\alpha _{i},\alpha _{j})>0$) and hence $G$
is a holomorphic elliptic function and is therefore constant in $z_{i}$.
Thus $G(z_{1},\ldots ,z_{n})=F_{N}(q)$ and the result follows. $\qed 
$

\end{document}